\documentclass[10pt]{amsart}

\usepackage{latexsym}
\usepackage{amsmath}
\usepackage{amssymb}
\begin{document}
 
%%%%%%%%%%%%%%%%%%%%%%%%%%%%%%%%%%%%%%%%%%%%%%%%%%%%%%%%%%%%
 
{\theoremstyle{plain}%
  \newtheorem{theorem}{Theorem}[section]
  \newtheorem{corollary}[theorem]{Corollary}
  \newtheorem{proposition}[theorem]{Proposition}
  \newtheorem{lemma}[theorem]{Lemma}
  \newtheorem{question}[theorem]{Question}
}
{\theoremstyle{remark}
\newtheorem{fact}{Fact}
\newtheorem{remark}[theorem]{Remark}
}

{\theoremstyle{definition}
\newtheorem{definition}[theorem]{Definition}
\newtheorem{example}[theorem]{Example}
}

% TOP MATTER

\title[The defining ideal of points in multi-projective space]% 
{On the defining ideal of a set of points in
multi-projective space} % This is the full title of the paper
% Avoid equations in title, but if you insist: $E=\lowercase{mc}^2$
% Do not use the \thanks{} command; use \extraline{} instead (see above).

\author{Adam Van Tuyl}
\thanks{Supported in part by the Natural Sciences and Engineering Research
Council of Canada}
\address{Department of Mathematical Sciences \\
Lakehead University \\
Thunder Bay, ON P7B 5E1, Canada}
\email{avantuyl@sleet.lakeheadu.ca}

%Insert `2000 Mathematics Subject Classification' numbers here!
\keywords{points, multi-projective space, generators}
\subjclass[2000]{13D02,13D40, 14A05}

\maketitle

% MACROS

\newcommand{\bj}{\underline{j}}
\newcommand{\bi}{\underline{i}}
\newcommand{\bt}{\underline{t}}
\newcommand{\dep}{\operatorname{depth}}
\newcommand{\hit}{\operatorname{ht}}
\newcommand{\op}{\overline{P}}
\newcommand{\xab}{X^{\underline{a}}Y^{\underline{b}}}
\newcommand{\Iz}{I_{\Z}}
\newcommand{\Ixp}{I_{\xp}}
\newcommand{\Z}{\mathbb{Z}}
\newcommand{\ax}{\alpha_{\X}}
\newcommand{\bx}{\beta_{\X}}
\newcommand{\kxo}{k[x_0,\ldots,x_n]}
\newcommand{\kx}{k[x_1,\ldots,x_n]}
\newcommand{\popo}{\mathbb{P}^1 \times \mathbb{P}^1}
\newcommand{\pr}{\mathbb{P}}
\newcommand{\pn}{\mathbb{P}^n}
\newcommand{\pnpm}{\mathbb{P}^n \times \mathbb{P}^m}
\newcommand{\pnk}{\mathbb{P}^{n_1} \times \cdots \times \mathbb{P}^{n_k}}
\newcommand{\X}{\mathbb{X}}
\newcommand{\Y}{\mathbb{Y}}
\newcommand{\N}{\mathbb{N}}
\newcommand{\M}{\mathbb{M}}
\newcommand{\Q}{\mathbb{Q}}
\newcommand{\Ix}{I_{\X}}
\newcommand{\pix}{\pi_1(\X)}
\newcommand{\pixt}{\pi_2(\X)}
\newcommand{\pipi}{\pi_1^{-1}(P_i)}
\newcommand{\piqi}{\pi_2^{-1}(Q_i)}
\newcommand{\qpi}{Q_{P_i}}
\newcommand{\pqi}{P_{Q_i}}
\newcommand{\pitk}{\pi_{2,\ldots,k}}
\newcommand{\kdim}{\operatorname{K-}\dim}
\newcommand{\reg}{\operatorname{reg}}
\newcommand{\regI}{\operatorname{reg({\it I}_\X)}}
\newcommand{\Hx}{\mathcal{H}_{\X}}
\newcommand{\A}{\mathcal{A}}
\newcommand{\D}{\mathcal{D}}
\newcommand{\E}{\mathcal{E}}
\newcommand{\T}{\mathcal{T}}
\newcommand{\B}{\mathcal{B}}
\newcommand{\sss}{\mathcal{S}}
\newcommand{\DD}{\mathbb{D}}
\newcommand{\F}{\mathcal{F}}

\begin{abstract}
We investigate the defining ideal $\Ix$ of a set of points $\X$ in $\pnk$
with a special emphasis on the case that $\X$ is in generic position,
that is, $\X$ has the maximal Hilbert function. 
When $\X$ is  in generic position,
we determine the degrees
of the  generators of the associated ideal $\Ix$.
Letting $\nu(\Ix)$ denote the minimal number of generators of $\Ix$,
we use this description of the degrees to construct a function
$v(s;n_1,\ldots,n_k)$ with the property that $\nu(\Ix) \geq 
v(s;n_1,\ldots,n_k)$ always holds for $s$ points in generic position
in $\pnk$.  When $k=1$, $v(s;n_1)$ equals the expected
value for $\nu(\Ix)$ as predicted by the Ideal
Generation Conjecture.  If $k\geq 2$, we show that
there are cases with $\nu(\Ix) > v(s;n_1,\ldots,n_k)$.  However,
computational evidence suggests that in many cases
$\nu(\Ix) = v(s;n_1,\ldots,n_k)$.
\end{abstract}

%%%%%%%%%%%%%%%%%%%%%%%%%%%%%%%%%%%%%%%%%%%%%%%%%%%%%%%%%%%%%%%%%%%%%%%

\section*{Introduction}

In this paper we investigate the generators of the ideal $\Ix$
defining a set of points $\X$ in generic position in $\pnk$.

One of the fundamental open problems about finite sets of 
points $\X \subseteq \pr^n$ in generic position, i.e.,
those sets of points having the maximal Hilbert function, is to
count the  minimal number of generators of 
$\Ix$ in terms of the data $n$ and $|\X| = s$.
This question is the content of
the Ideal Generation Conjecture (IGC) (see \cite{GO2}). 
Recently, many authors (cf. 
\cite{CGG,GuMaRa1,GuMaRa2,GuMaRa3,Gu,GVT,VT1,VT2}) 
have been interested in generalizing
results about points in $\pr^n$ to $\pnk$.  
We continue this program by studying the
generators of $\Ix$ when $\X \subseteq \pnk$
with the hope that this might lead to
a generalized IGC. 
Our investigation was also partially motivated
by the desire to understand which properties about
the ideal of points in $\pr^n$,
specifically those shown in \cite{GM,GO,GO2}, carry over to $\pnk$.

Given the defining ideal $\Ix$ of a set of points $\X \subseteq \pnk$,
two natural questions about the generators of $\Ix$ arise: $(1)$
what are the degrees of the generators? and $(2)$ what is
$\nu(\Ix) :=$ minimal number of generators of $\Ix$?  These
questions can be viewed as the first step in describing the 
multi-graded minimal free resolution of $\Ix$ since $(1)$
and $(2)$ are questions about the $0$th multi-graded
Betti numbers.

In Section 2 we show that the Hilbert function of a set of points
can be used to bound the degrees of the generators, thus giving a 
partial answer to $(1)$.  As posed, however, these questions
are difficult to attack, even when $k=1$, without further conditions on
the points.

For finite sets of points $\X \subseteq \pr^n$, these questions
have been primarily studied under the extra hypothesis that the
set of points is in generic position, i.e.,
$H_{\X}(i) = \min \{\dim_{\bf k} R_i, |\X|\}$ for all $i \in \N$.
Thus, one is led to ask about the generators of $\Ix$ when
$\X$ is a set of points in generic position in $\pnk$.  However,
it is first necessary to establish the basic properties (like existence)
of points in generic position in multi-projective spaces since these
facts are not part of the literature.  Analogous to the case
of points in $\pr^n$, we say that a set $\X$ of $s$ points
in $\pnk$ is in {\it generic position} if $H_{\X}(\bi) = 
\min\{\dim_{\bf k} R_{\bi}, s\}$ for all $\bi \in \N^k$.  In Section 3
we show that these points exist, and moreover, if we consider each
set of $s$ points as a point in $(\pnk)^s$, the points
in generic position form a non-empty open subset of $(\pnk)^s$
with respect to the Zariski topology.

We also show in Section 3 that if $F$ is a generator of $\Ix$,
then $\deg F = \bi$ or $\bi + e_j$ where $\bi \in \D :=
\min\left\{\bi \in \N^k 
~\left|~ \binom{i_1+n_1}{i_1}\cdots\binom{i_k+n_k}{i_k} >s 
\right\}\right.$ 
and $e_j$ is one of the $k$ basis vectors of $\N^k$.  This result
gives an answer to $(1)$ and generalizes the fact that $\Ix =
\langle (\Ix)_d \oplus (\Ix)_{d+1} \rangle$ with 
$d = \min\left\{ i ~\left|~ \binom {i+n}{i} > s \right\}\right.$ in
the graded case.  An interesting difference between points
in generic position in $\pr^n$ versus $\pnk$
is that $R/\Ix$ is {\it always} Cohen-Macaulay if $k=1$, but 
is {\it never} Cohen-Macaulay if $k \geq 2$ (see Theorem 
\ref{depthofgenericpoints})

In Section 4 we use this description of the degrees to show that
$\nu(\Ix)$ can be determined by counting the generators of 
degree $\bi$ and $\bi + e_j$ for all $\bi \in \D$ and 
$j = 1,\ldots,k$.  By degree considerations, $\Ix$ has
$\dim_{\bf k} (\Ix)_{\bi}$ generators of degree $\bi$
for each $\bi \in \D$.  To count the generators of degree
$\bi + e_j$, we need to calculate the dimension
of the image of the map $\Phi_{\bi,j} : R_{e_j} \otimes_{\bf k}
(\Ix)_{\bi} \stackrel{a \times b}{\longrightarrow} (\Ix)_{\bi + e_j}$
for each $\bi \in \D$ and $1 \leq j \leq k$.  Moreover, if there
exists $\bi_1,\bi_2 \in \D$ and $1 \leq j_1,j_2 \leq k$
such that $l = \bi_1 + e_{j_1} = \bi_2 + e_{j_2}$, then
Im $\Phi_{\bi_1,e_{j_1}}$ and Im $\Phi_{\bi_1,e_{j_2}}$ are
both subspaces of $(\Ix)_l$, so we also need to know
$\dim_{\bf k}(\operatorname{Im} _{\bi_1,e_{j_1}} \cap \operatorname{Im} 
\Phi_{\bi_1,e_{j_2}})$. In general, it is difficult to compute
the sizes of these vector spaces, even if $k=1$, except in some 
special cases.  When $k=1$, to compute $\nu(\Ix)$ only the dimension
of the image of $\Phi : R_1 \otimes_{\bf k} (\Ix)_d 
\rightarrow (\Ix)_{d+1}$ needs to be calculated.  The IGC states that
Im $\Phi$ should be as large as possible for a sufficiently
general set of points (a subset of those points in generic position).

By considering the largest possible value for each $\dim_{\bf k} 
\operatorname{Im} \Phi_{\bi,j}$, in Section 5 we construct
a function $v(s;n_1,\ldots,n_k)$ with the
property that $\nu(\Ix) \geq v(s;n_1,\ldots,n_k)$ always holds
for a set of $s$ points in generic position in $\pnk$.  When
$k=1$, $v(s;n)$ equals the expected value for $\nu(\Ix)$ as
predicted by the IGC.  Furthermore, using \cite{GuMaRa2,GuMaRa3},
we show that $\nu(\Ix) = v(s;1,1)$ for a sufficiently general
set of $s$ points in $\popo$.

Buoyed by these results, we had hoped that for any set of $s$ points
in generic position that were sufficiently general, we should
expect $\nu(\Ix) = v(s;n_1,\ldots,n_k)$.  However, we show
that if $\X$ is any three points in generic position in 
$\pr^1 \times \cdots \times \pr^1$ ($k\geq 3$ times), then
$\nu(\Ix) > v(3;1,\ldots,1)$.  As well, computational
evidence suggests that $\nu(\Ix) > v(s;1,n,n)$ if $s = 1+n+n$.  
These cases appear to be exceptional because in all other computed
examples the equality $\nu(\Ix) = v(s;n_1,\ldots,n_k)$ holds.
Moreover, we know of no example of $\nu(\Ix) > v(s;n_1,n_2)$
when $k=2$.  This leads us to believe
that $\nu(\Ix) = v(s;n_1,\ldots,n_k)$ in a large number of cases, thus
giving us a partial generalization of the the IGC.

%%%%%%%%%%%%%%%%%%%%%%%%%%%%%%%%%%%%%%%%%%%%%%%%%%%%%%%%%%%%%%%%%%%%%%%

\section{Preliminaries}

In this paper {\bf k} denotes a field with char$({\bf k}) = 0$
and ${\bf k} =\overline{\bf k}$.
This section provides the relevant facts and definitions
about multi-graded rings, Hilbert functions,  and sets
of points in multi-projective spaces.
See also \cite{Ro,VT1,VT2}.

Let $\N:=\{0,1,2,\ldots\}$.  For any integer $k \geq 1$,
we write $[k]:=\{1,\ldots,k\}$.  We denote $(i_1,\ldots,i_k) \in \N^k$
by $\bi$.  
We set $|\bi| := \sum_h i_h$.
If $\bi,\bj \in \N^k$, then $\bi + \bj := (i_1 + j_1, \ldots,i_k + j_k)$.
We write $\bi \geq \bj$ if $i_h \geq j_h$ for every $h= 1,\ldots,k$. 
Observe that $\geq$ is a partial order on $\N^k$.  For any subset
$\A \subseteq \N^k$, we will use $\min \A$ to denote the
set of minimal elements of $\A$ with respect to this
partial order. 
The set $\N^k$ is a semi-group generated by $\{e_1,\ldots,e_k\}$
where $e_i:= (0,\ldots,1\ldots,0)$ is the $i$th 
standard basis vector of $\N^k$.  For any $c \in \N$, 
$ce_i := (0,\ldots,c,\ldots,0)$.

Set  $R = {\bf k}[x_{1,0},\ldots,x_{1,n_1},x_{2,0},\ldots,x_{2,n_2},
\ldots,x_{k,0},\ldots,x_{k,n_k}]$, and induce an $\N^k$-grading on $R$
by setting $\deg x_{i,j} = e_i$.  An element $x \in R$ is said to
be $\N^k$-{\it homogeneous} (or simply {\it homogeneous}
if the grading is clear) if $x \in R_{\bi}$
for some $\bi \in \N^k$.  If $x$ is homogeneous, then
$\deg x := \bi$.

An ideal $I = (F_1,\ldots,F_r) \subseteq R$ is an  $\N^k$-{\it homogeneous}
(or, simply {\it homogeneous}) {\it ideal} if each
$F_j$ is $\N^k$-homogeneous.  If $I \subseteq R$ is a homogeneous ideal,
$S = R/I$ inherits an
$\N^k$-graded ring structure if we define $S_{\bi} = (R/I)_{\bi} :=
R_{\bi}/I_{\bi}$.   The function
$H_S(\bi) := \dim_{\bf k} (R/I)_{\bi}$ 
is the {\it Hilbert function of $S$}.  

If $I$ is an $\N^k$-homogeneous ideal of $R$, then for
any $\bi \in \N^k$, and for any $j \in [k]$, we set
\[
R_{e_j}I_{\bi} := \left\{f ~\left|~ f = f_0x_{j,0} + f_1x_{j,1} +
\cdots + f_{n_j}x_{j,n_j}, ~f_l \in I_{\bi}\right\}\right..
\]
Note that $R_{e_j}I_{\bi}$ is a subspace of the vector
space $I_{\bi+e_j}$.

For every $\bi \in \N^k$, a basis for $R_{\bi}$ as a vector space 
over {\bf k} is the set of 
all monomials in $R$ of degree $\bi$.
Thus, $\dim_{\bf k} R_{\bi} = 
\binom{n_1+i_1}{i_1}\binom{n_2 + i_2}{i_2}\cdots \binom{n_k + i_k}{i_k}$.
We set $N(\bi) := \dim_{\bf k} R_{\bi}$ for each $\bi \in \N^k$.
 
The $\N^k$-graded ring $R$ 
is the coordinate ring of $\pnk$. If 
$P \in \pnk$
is a point, and  if $I_P$ denotes the ideal associated to  $P$,  
then the ideal $I_P$ is a prime ideal, and furthermore,
$
I_P = (L_{1,1},\ldots,L_{1,n_1},\ldots,L_{k,1},\ldots,L_{k,n_k})
$
where $\deg L_{i,j} = e_i$ for $j = 1,\ldots,n_i$.
Let  $P_1,\ldots,P_s$ be $s$ distinct points 
in $\pnk$.  If $\X = \{P_1,\ldots,P_s\}$, then 
the $\N^k$-homogeneous ideal $\Ix$ of
forms that vanish at $\X$ is $\Ix = I_{P_1} \cap \cdots \cap I_{P_s}$ where
$I_{P_i}$ is the ideal associated to the point $P_i$.  
The coordinate
ring  $R/\Ix$ then has the following property.

\begin{lemma}[{\cite[Lemma 3.3]{VT1}}] \label{firstnzd1}
If $\X \subseteq \pnk$ is a finite set  of points, then   
for each $l \in [k],$  there exists $L_l \in R_{e_l}$  
such that $\overline{L}_l$  is a non-zero divisor in $R/\Ix$.
\end{lemma}

\begin{remark} \label{nzd}
After a linear change of variables in the $x_{1,j}$'s, a change
of variables in the $x_{2,j}$'s, and so on, we can take
$L_l = x_{l,0}$ for each $l \in [k]$.  We 
therefore assume, once and for all,
that the set of points $\X$ under investigation has
the property that $\overline{x}_{l,0}$ is a non-zero divisor in $R/\Ix$
for each $l \in [k]$.
\end{remark} 

We sometimes write $H_{\X}$ for $H_{R/\Ix}$,
and  call $H_{\X}$ the Hilbert function of $\X$. 
Classifying the Hilbert functions   
of sets of points in $\pnk$ with $k \geq 2$ remains an open problem
(the case $k=1$ is dealt with in \cite{GeMaRo}).  See
\cite{GuMaRa1} and \cite{VT1} for some progress on this problem.  
However, some growth conditions on $H_{\X}$ are known.

\begin{theorem}[{\cite[Proposition 3.5]{VT1}}]
\label{simpleboundsonhx}
Let $\X$ be a finite set of points in $\pnk$ with Hilbert
function $H_{\X}$.
\begin{enumerate}
\item[$(i)$] For all $\bi \in \N^k$, 
$H_{\X}(\bi) \leq H_{\X}(\bi + e_l)$ for all $l \in [k]$. 
\item[$(ii)$] If $H_{\X}(\bi)
=H_{\X}(\bi + e_l)$  for some $l \in [k]$, then 
$H_{\X}(\bi + e_l)
=H_{\X}(\bi + 2e_l).$
\end{enumerate}
\end{theorem}

Let $\pi_i:\pnk \rightarrow \pr^{n_i}$ denote the $i$th
projection morphism defined by $P_1 \times \cdots \times P_i \times
\cdots \times P_k \mapsto P_i$.  Then $\pi_i(\X)$ is
the set of  all the $i$th-coordinates in $\X$.  For each $i \in [k]$,
 set $t_i := |\pi_i(\X)|$. With this
notation we have

\begin{theorem}[{\cite[Corollary 4.7]{VT1}}] \label{eventualgrowth}
Let $\X$ be a finite set of  points
in $\pnk$ with Hilbert function $H_{\X}$.   Fix an $i \in [k]$.
If $(j_1,\ldots,j_i,\ldots,j_k) \in \N^k$ with $j_i \geq t_i - 1$,
then 
$ 
H_{\X}(j_1,\ldots,j_i,\ldots,j_k) = H_{\X}(j_1,\ldots,t_i-1,\ldots,j_k).
$
\end{theorem}

\begin{remark}One can interpret the above results as follows.  
Fix an integer $i \in [k]$, and fix $k-1$ integers in $\N$,
say $j_1,\ldots,j_{i-1},j_{i+1},\ldots,j_k$.  Set
\[\underline{j}_{l} := (j_1,\ldots,j_{i-1},l,j_{i+1},\ldots,j_k)
~~\mbox{for each integer $l \in \N$.}\]  
Then Theorems \ref{simpleboundsonhx}
and \ref{eventualgrowth} imply that
there exists an integer $l' \leq t_i - 1$ such that the
sequence $
H_{\X}(\underline{j}_0),H_{\X}(\underline{j}_1),
H_{\X}(\underline{j}_2),
H_{\X}(\underline{j}_3),....$
has the property that $H_{\X}(\underline{j}_{l})
< H_{\X}(\underline{j}_{l+1})$ if $0 \leq l < l'$, but
$H_{\X}(\bj_{l}) = H_{\X}(\bj_{l+1})$ if $l \geq l'$.
\end{remark}

%%%%%%%%%%%%%%%%%%%%%%%%%%%%%%%%%%%%%%%%%%%%%%%%%%%%%%%%%%%%%%%%%%%%%%%

\section{On the generators of an ideal of a set of points}

Let $\Ix$ be the defining ideal of a finite set of points $\X
\subseteq \pnk$ with Hilbert function $H_{\X}$.   Using only 
 $H_{\X}$,
we describe a finite subset $\E \subseteq \N^k$ with the property
that if $F$ is a 
generator of $\Ix$, then $\deg F \in \E$.  

Fix an $l \in [k]$.  Then, for each $\bj=(j_1,\ldots,j_{l-1},j_{l+1},
\ldots,j_k) \in \N^{k-1}$, set
\[
i(\bj) := \min 
\left\{ i \in \N^+ ~\left|~ 
\begin{array}{c}
H_{\X}(j_1,\ldots,j_{l-1},i-1,j_{l+1},\ldots,j_k) =  \\
H_{\X}(j_1,\ldots,j_{l-1},i,j_{l+1},\ldots,j_k)
\end{array}\right\}\right..
\]
The existence of the integer $i(\bj)$ follows from
Theorem \ref{eventualgrowth}. 

\begin{theorem}	\label{generators}\label{nomingen1}
Let $\X$ be a finite set of points
of  $\pnk$.  Fix an
$l \in [k]$ and $\bj = (j_1,\ldots,j_{l-1},j_{l+1},
\ldots,j_k) \in \N^{k-1}$.
Set $\bi = (j_1,\ldots,j_{l-1},i(\bj),j_{l+1},\ldots,j_k)$.
Then
\[
(\Ix)_{\bi+(r+1)e_l}  = 
R_{e_l} (\Ix)_{\bi+re_l} ~~\mbox{ for all $r \in \N$.} \]
In particular, 
if there exists $\underline{l} \in \N^k$ and $t \in [k]$
such that $H_{\X}(\underline{l}) = H_{\X}(\underline{l} - e_t) = 
H_{\X}(\underline{l}-2e_t)$,
then $\Ix$ has no minimal generators of degree $\underline{l}$.
\end{theorem}

\begin{proof}  
Without loss of generality, we only consider 
the case $l = 1$.  By Remark \ref{nzd} we can take
$x_{1,0}$ to be a non-zero divisor.  Set 
$S = {\bf k}[x_{1,1},\ldots,x_{k,n_k}]
\cong R/(x_{1,0})$ and $\bi = (i(\bj),j_2,\ldots,j_k)$ where
$\bj = (j_2,\ldots,j_k)$.  

Because $x_{1,0}$ is a non-zero divisor, $\Ix/x_{1,0}\Ix 
\cong \left((\Ix,x_{1,0})/x_{1,0}\right)$.  
For each $\bt \in \N^k$ the short exact sequence
\[ 
0 \longrightarrow (\Ix)_{\bt -e_1} \stackrel{\times x_{1,0}}
{\longrightarrow}
(\Ix)_{\bt} \longrightarrow \left(\Ix/(x_{1,0}\Ix)\right)_{\bt}
\longrightarrow 0 
\]
implies $\dim_{\bf k} (\Ix)_{\bt} = \dim_{\bf k} (\Ix)_{\bt -e_1}
+ \dim_{\bf k} \left((\Ix,x_{1,0})/x_{1,0}\right)_{\bt}$.
On the other hand, the short exact sequence
\[0 \longrightarrow R/\Ix(-e_1)  \stackrel{\times \overline{x}_{1,0}}
{\longrightarrow} R/\Ix \longrightarrow R/(\Ix,x_{1,0})
\cong \frac{R/(x_{1,0})}{(\Ix,x_{1,0})/x_{1,0}} \longrightarrow 0,\]
and the hypothesis that $\dim_{\bf k} (R/\Ix)_{\bi + re_1}
= \dim_{\bf k} (R/\Ix)_{\bi + (r-1)e_1}$ for every $r \in \N$
implies that $\left((\Ix,x_{1,0})/x_{1,0}\right)_{\bi +re_1} 
= (R/(x_{1,0}))_{\bi +re_1} \cong
S_{\bi +re_1}$ for all $r \in \N$.

Fix an integer $r \in \N$ and set $W = R_{e_1}(\Ix)_{\bi+re_1}$.
Then $W \subseteq (\Ix)_{\bi + (r+1)e_1}$,
and because $x_{1,0}$ is a non-zero divisor
\begin{eqnarray*}
\dim_{\bf k} W & = &\dim_{\bf k} (\Ix)_{\bi+ re_1} + \dim_{\bf k} W'\\
&=& \dim_{\bf k} (\Ix)_{\bi + (r+1)e_1} - \dim_{\bf k} S_{\bi + (r+1)e_1}
+ \dim_{\bf k} W'
\end{eqnarray*}
where $W' = \{f'' ~|~ f = f'x_{1,0} + f'', ~f \in W\}
= \{f(0,x_{1,1},\ldots,x_{k,n_k}) ~|~ f \in W\}$.  By slightly abusing
notation,
$W'$ can be viewed as a subset of $S_{\bi+(r+1)e_1}$.

It suffices to show that $S_{\bi + (r+1)e_1} \subseteq W'$
because then $\dim_{\bf k} S_{\bi + (r+1)e_1}
= \dim_{\bf k} W'$, 
whence $\dim_{\bf k} W = \dim_{\bf k} (\Ix)_{\bi +(r+1)e_1}$
which gives $W = (\Ix)_{\bi + (r+1)e_1}$.  So,
suppose $f \in S_{\bi + (r+1)e_1}$.  Then
$f = f_1x_{1,1} + \cdots + f_{n_1}x_{1,n_1}$ with $f_i \in S_{\bi + re_1}$.
But $S_{\bi + re_1} \cong 
\left((\Ix,x_{1,0})/{x_{1,0}}\right)_{\bi + re_1},$
so by abusing notation, there exists $F_i \in (\Ix)_{\bi+re_1}$
such that $F_i = g_ix_{1,0} + f_i$.  But then 
$F = F_1x_{1,1} + \cdots +F_{n_1}x_{1,n_1} \in W$, whence $f \in W'$.

For the last statement let $\underline{l} = (l_1,\ldots,l_k)$
and  $i:=i(l_1,\ldots,l_{t-1},l_{t+1},\ldots,l_k)$.
Then $H_{\X}(\underline{l}) = H_{\X}(\underline{l}-e_t) =
H_{\X}(\underline{l}-2e_t)$ implies that $\underline{l} - e_t 
\geq (l_1,\ldots.l_{t-1},i,l_{t+1},\ldots,l_k)$,
so $(\Ix)_{\underline{l}} = R_{e_t}(\Ix)_{\underline{l} - e_t}$.
\end{proof}

Let $\X$ be a finite
set of points in $\pnk$, and set 
\[\B := 
\left\{\bi\in\N^k ~\left|~ H_{\X}(\bi) < \dim_{\bf k} R_{\bi} = N(\bi) \right\}
\right.
= \left\{\bi\in\N^k ~\left|~ (\Ix)_{\bi} \neq 0 \right\}\right..
\]
For each $l \in [k]$ and for each  
$\bj = (j_1,\ldots,j_{l-1},j_{l+1},\ldots,j_k) \in \N^{k-1}$  set 
\[
\A_{l,\bj} := \left\{(j_1,\ldots,j_{l-1},i(\bj),j_{l+1},\ldots,j_k)
+ re_l ~\left|~ r \in \N^+ \right\}\right..
\]
We then define 
\[\A := \bigcup_{l=1}^k\left(\bigcup_{\bj\in\N^{k-1}} \A_{l,\bj}\right).\]

\begin{theorem} \label{generatorsprop}
Let $\X$ be a finite set of  points in $\pnk$ with defining ideal $\Ix$.
With the notation as above, set $\E = \B \backslash \A$.
Then $\E$ is a finite set.  Furthermore, if $f$ is a  generator
of $\Ix$, then $\deg f  \in \E$.  In particular,
$\Ix = \left\langle \bigoplus_{\bi \in \E} (\Ix)_{\bi} \right\rangle.$
\end{theorem}

\begin{proof}
We show that $\E$ is finite.  Let $t_i := |\pi_i(\X)|$ where 
$\pi_i:\pnk \rightarrow \pr^{n_i}$ is the $i$th projection morphism.
Set $\F := \{\bi \in \N^k ~|~ \bi \leq (t_1,t_2,\ldots,t_k)\}$.
Now suppose that $\bj = (j_1,\ldots,j_k) \in \N^k \backslash \F$.
Thus, there is a coordinate of $\bj$, say $j_i$, such that
$j_i \geq t_i + 1$.  By Theorem \ref{eventualgrowth}
\[
H_{\X}(\bj) = H_{\X}(j_1,\ldots,j_{i-1},t_i,j_{i+1},\ldots,j_k)
= H_{\X}(j_1,\ldots,j_{i-1},t_i-1,j_{i+1},\ldots,j_k).
\]
This means that $\bj \in \A_{i,(j_1,\ldots,j_{i-1},j_{i+1},\ldots,j_k)}
\subseteq \A$.  We thus have $\N^k \backslash \F \subseteq \A$. 
But this implies that $\E = \B \backslash \A \subseteq 
\N^k \backslash \A \subseteq \F$, and since $\F$ is finite,
so is $\E$.

For the second statement, let $f$ be a generator of $\Ix$.
Then it is immediate that $\deg f \in \B$.  On the other hand,
Theorem \ref{generators} implies that $\deg f \not\in \A_{l,\bj}$ for
any $l \in [k]$ or $\bj \in \N^{k-1}$.  Hence,
$\deg f \not\in \A$, so $\deg f \in \E$.
\end{proof}

\begin{remark} 
We recover Proposition 1.1 (3) of \cite{GM}
when $k=1$, i.e., $\X \subseteq \pr^n$.
\end{remark}

%%%%%%%%%%%%%%%%%%%%%%%%%%%%%%%%%%%%%%%%%%%%%%%%%%%%%%%%%%%%%%%%%%%%%%%

\section{Points in generic position}

Analogous to the definition for points in $\pr^n$, a 
set of points in $\pnk$ is said to be in generic position 
if its Hilbert function is maximal.  Although such sets have
been studied (cf. \cite{GuMaRa2,GuMaRa3}) we could find
no proof in the literature for the existence of such sets
when $k \geq 2$ (the case $k = 1$ is \cite[Theorem 4]{GO}).
We therefore begin by providing a proof of this ``folklore'' result.
Then, if $\Ix$ is the defining ideal of a set of points in
generic position in $\pnk$, we compute the depth of $R/\Ix$, and 
give bounds on the degrees of the generators of $\Ix$.

Theorems \ref{simpleboundsonhx}
and \ref{eventualgrowth} imply the number of possible Hilbert functions 
for $s$ points is finite.
However, since the number of sets with $s$ points is infinite, 
we can ask if there exists an expected Hilbert function for
$s$ points.   We give a heuristic argument for this expected
function.

If $\{m_1,\ldots,m_{N(\bj)}\}$ are the $N(\bj)$ monomials of degree
$\bj$ in the $\N^k$-graded ring $R$, then any  
$F \in R_{\bj}$ can be written as 
$F = \sum_{i=1}^{N(\bj)} c_im_i$ where $c_i \in {\bf k}$.
Suppose that $P \in \pnk$.  For
$F \in R_{\bj}$ to vanish at $P$ we require
$
F(P) = \sum_{i=1}^{N(\bj)} c_im_i(P) = 0.
$
By considering the $c_i$'s as unknowns, this equation
gives us one linear condition.
If $\X = \{P_1,\ldots,P_s\}$, then for $F \in R_{\bj}$
to vanish on $\X$ we  require that $F(P_1) = \cdots
= F(P_s) = 0$.  We then have a linear system of equations
\[\bmatrix
m_1(P_1) & \cdots & m_{N(\bj)}(P_1) \\
\vdots &  &\vdots \\
m_1(P_s) & \cdots & m_{N(\bj)}(P_s) 
\endbmatrix
\bmatrix
c_1 \\
\vdots \\
c_{N(\bj)} 
\endbmatrix
=
\bmatrix
0 \\
\vdots \\
0
\endbmatrix.
\]
The number of linearly independent solutions is the rank of the matrix
on the left.  For a general enough set of points, we expect
this rank to be as large as possible.  By \cite[Proposition 4.3]{VT1} 
the rank of this matrix equals  $H_{\X}(\bj)$,
so we expect a general enough set of $s$ points 
$\X \subseteq \pnk$ to have the Hilbert function
$H_{\X}(\bj) = \min\left\{ N(\bj),s \right\}$ for all
$\bj \in \N^k.$
Proceeding analogously as in the case of points in $\pr^n$, we make
the following definition.

\begin{definition}
Let $\X$ be a finite set of $s$ points in  $\pnk$ 
with Hilbert function $H_{\X}$.  If
\[H_{\X}(\bj) = \min\left\{ N(\bj),s \right\} 
\hspace{.5cm}\text{for all}~
\bj \in \N^k,\]
then the Hilbert function is called {\it maximal}.  
A set of $s$ points is said to be in {\it generic position}
if its Hilbert function is maximal.
\end{definition}

We now show the existence of 
sets of points in generic position by demonstrating that ``most'' sets
of $s$ points in $\pnk$ are in generic position.  
We shall denote $(\pnk)\times \cdots \times
(\pnk)$ ($s$ times) by $(\pnk)^s$.

\begin{theorem}
\label{genericpoints}
The $s$-tuples of points of $\pnk$, $(P_1,\ldots,P_s)$, considered as points
of $(\pnk)^s$, which are in generic position form a non-empty open
subset of $(\pnk)^s$.
\end{theorem}

\begin{proof}
Since the case $k=1$ is found in \cite{GO}, we can assume that $k \geq 2$.
Let $\{m_1,\ldots,m_{N(\bj)}\}$ be the $N(\bj)$ monomials
of degree $\bj = (j_1,\ldots,j_k) \in \N^k$ in $R$.  By composing the
product of $j_i$-uple embeddings with the Segre embedding we have a morphism
$\nu_{\bj}:\pnk \rightarrow \pr^{N(\bj)-1}$ 
defined by
\[ P = [a_{1,0}:\cdots:a_{1,n_1}]\times \cdots \times
 [a_{k,0}:\cdots:a_{k,n_k}] \longmapsto [m_1(P):\cdots:m_{N(\bj)}(P)],\]
i.e., $m_i(P)$ is the monomial $m_i$ evaluated at $P$.
This induces a morphism
\[\varphi_{\bj} = (\nu_{\bj}:\cdots:\nu_{\bj}):\left(\pnk\right)^s
\longrightarrow \left(\pr^{N(\bj)-1}\right)^s = V_{\bj}.\]
By \cite{GO} there exists a nonempty open subset $W_{\bj} \subseteq
V_{\bj}$ with the property that each point $(Q_1,\ldots,Q_s)
\in W_{\bj}$ corresponds to a set of $s$ points in $\pr^{N(\bj)-1}$
in generic position.  In particular, each point of $W_{\bj}$
corresponds to a set of $s$ points in $\pr^{N(\bj)-1}$ that
impose $\min\{s,N(\bj)\}$ conditions on linear forms.

Because $\nu_{\bj}$ does not vanish everywhere, $U_{\bj} :=
\varphi^{-1}(W_{\bj})$ is a non-empty open subset of $(\pnk)^s$.
Furthermore, because $\nu_{\bj}$ induces an isomorphism
between the linear forms of $\pr^{N(\bj)-1}$ and the forms
of $R_{\bj}$, if $(P_1,\ldots,P_s) \in U_{\bj}$, then
as a set of points $\X = \{P_1,\ldots,P_s\}$ of
$\pnk$ we have $H_{\X}(\bj) = \min\{s, N(\bj)\}$.  Hence
$U = \bigcap_{\bj \in \N^k} U_{\bj}$ consists of those $s$-tuples which 
correspond to sets of $s$ points in $\pnk$ with maximal Hilbert functions.

To complete the proof, it suffices to show that the intersection
$U = \bigcap_{\bj \in \N^k} U_{\bj}$ can be taken to be finite, and thus
$U$ is open.  Suppose $\bj \in \N^k$ is such that $s = \min\{s,N(\bj)\}$.
We claim that $U_{\bj} \subseteq U_{\bj'}$ for all $\bj \leq \bj'$.
Indeed, take $(P_1,\ldots,P_s) \in U_{\bj}$.  So, if $\X =
\{P_1,\ldots,P_s\}$, we have $H_{\X}(\bj) = s$.  Since the
Hilbert function strictly increases until it stabilizes
and is bounded by $s$, $H_{\X}(\bj') = s < N(\bj')$ for
all $\bj \leq \bj'$.  Thus $(P_1,\ldots,P_s) \in U_{\bj'}$
as desired.  Setting $D_1 = \{\bj\in\N^k ~|~ N(\bj) < s\}$
and $D_2 = \min\{\bj \in \N^k ~|~ N(\bj) \geq s\}$,
we thus have $D = D_1 \cup D_2$ is a finite set
and $U = \bigcap_{\bj \in D} U_{\bj}$.
\end{proof}

For any finite set of points $\X \subseteq \pnk$,
we have $\operatorname{K-}\dim 
R/\Ix = k$.  However, it was shown in \cite[Proposition 2.6]{VT2}
that the depth of $R/\Ix$ may take on any value in $\{1,\ldots,k\}$.
When $\X$ is in generic position the depth can be calculated.
We begin with a lemma.

\begin{lemma}
\label{binomid}
Let $n,l \geq 1$ be integers.  Then 
$\binom{n+l+1}{l+1} \leq \binom{n+l}{l}(n+1)$.
\end{lemma}

\begin{proof}
Note that 
$\binom{n+l+1}{l+1} = \binom{n+l}{l}\cdot\frac{(n+l+1)}{(l+1)} 
= \binom{n+l}{l}\left(1 + \frac{n}{l+1}\right).$ The
inequality now follows since 
$l \geq 1$,  and thus $(1 + \frac{n}{l+1}) \leq (1+n)$.  
\end{proof}

\begin{theorem}	\label{depthofgenericpoints}
If $\X$ is a set of $s > 1$ points in generic position
in $\pnk$,
then $\operatorname{depth} R/\Ix = 1$.  In particular,
$R/\Ix$ is Cohen-Macaulay if and only if $k = 1$.
\end{theorem}

\begin{proof}
By Lemma \ref{firstnzd1}  $\operatorname{depth} R/\Ix \geq 1$.
We show that equality holds.  Without loss of generality, 
take $n_1 \leq n_2 \leq \cdots \leq n_k$ and  let $l$ be
the minimal integer such that $\binom{n_1+l}{l} \geq s$.  Since
$\X$ is in generic position,
$H_{\X}(l,0,\ldots,0) = \min\left\{\binom{n_1+l}{l},s\right\} =s.$

\noindent
{\it Claim. }  If $\bj \in \N^k$ and $\bj > (l-1,0,\ldots,0)$, then
$H_{\X}(\bj) =s$.

\noindent
{\it Proof of the Claim. }  
If $j_1 > l-1$, then $\binom{n_1+ j_1}{j_1} \geq \binom{n_1 + l}{l}$.
Thus 
$N(\bj) \geq \binom{n_1+l}{l}
\geq \min\left\{\binom{n_1+l}{l},s\right\} =s,$
and hence, $H_{\X}(\bj) =s$.
So, suppose $j_1 = l-1$.  Since $\bj > (l-1,0,\ldots,0)$, there
exists $m \in \{2,\ldots,k\}$ such that $j_m > 0$.  Since
$n_1 \leq n_m$, we have the following inequalities:
$N(\bj) \geq \binom{n_1+l-1}{l-1}\binom{n_m+j_m}{j_m} \geq 
\binom{n_1 + l-1}{l-1}\binom{n_1 + 1}{1}.$
By Lemma \ref{binomid}, we also have $\binom{n_1+l-1}{l-1}(n_1 +1)
\geq \binom{n_1 +l}{l}$.  Hence,
$N(\bj) \geq \binom{n_1+l}{l}
\geq \min\left\{\binom{n_1+l}{l},s\right\} =s.$
Therefore, $H_{\X}(\bj)=s$, as desired.\hfill$\Box$

Since $\overline{x}_{1,0}$ is a non-zero divisor of $R/\Ix$
we have the short exact sequence
\[
0 \longrightarrow \left(R/\Ix\right)(-e_1) 
\stackrel{\times \overline{x}_{1,0}}{\longrightarrow}
R/\Ix \longrightarrow R/(\Ix,x_{1,0}) = R/J \longrightarrow 0
\]
where $J = (\Ix,x_{1,0})$.  
Thus the Hilbert function of $R/J$ is
$H_{R/J}(\bj) = H_{\X}(\bj) - H_{\X}(\bj -e_1)$
for all $\bj \in \N^k$,
where $H_{\X}(\bj) = 0$ if $\bj \not\geq \underline{0}.$
From the claim, it follows that if $\bj > (l,0,\ldots,0)$,
then 
\[H_{R/J}(\bj) = H_{\X}(\bj) - H_{\X}(\bj -e_1)  = s - s = 0.\]
On the other hand, if $\bj = (l,0,\ldots,0)$, then
\[H_{R/J}(\bj) = H_{\X}(le_1) - H_{\X}((l-1)e_1) = 
s - \binom{n_1 + l-1}{l-1} > 0.\]
Since $s \neq 1$ there exists a non-constant element
 $F \in R_{le_1}$ such that 
$0 \neq \overline{F} \in R/J$.

It suffices to demonstrate that
all the non-constant homogeneous elements of $R/J$
are annihilated by $\overline{F}$, and hence,
$\operatorname{depth} R/J = 0$.  So, suppose
that $G \in R$ is such that $0 \neq \overline{G} \in R/J$.
Without loss of generality we can take $G$ to be an
$\N^k$-homogeneous element with $\deg G = (j_1,\ldots,j_k)> \underline{0}$.
Now $\deg FG = 
(j_1 + l,j_2,\ldots,j_k) > (l,0,\ldots,0)$.  Since
$H_{R/J}(j_1 + l,j_2,\ldots,j_k) = 0$, it follows that $FG \in J$.
\end{proof}

\begin{remark}
If $s=1$, then $\operatorname{depth} R/\Ix = k$
because the ideal of a point is a complete intersection.
\end{remark}

We now apply Theorem \ref{generatorsprop} to describe
the degrees of the generators of $\Ix$
when $\X$ is in generic position.
We introduce
some notation:  if $E = \{\bj_1,\ldots,\bj_l\} \subseteq \N^k$,
and $\bi \in \N^k$, then
\[
E + \bi := 
\left\{\bj_1 + \bi, \bj_2+\bi,\ldots,\bj_l + \bi\right\} \subseteq \N^k.
\]
Also, let $\D := \min\{\bi \in \N^k ~|~ (\Ix)_{\bi} \neq 0\} = 
\min\{\bi \in \N^k ~|~ H_{\X}(\bi) < N(\bi)\}.$  
(Recall that $\min \mathcal{S}$ with $\mathcal{S} \subseteq \N^k$
is the set of minimal elements of $\mathcal{S}$ with respect
to the partial ordering $\bi \geq \bj$ if $i_h \geq j_j$
for all $h$.)
When $\X$ is a set of $s$ points in generic position we have
$\D = \min\{\bi \in \N^k ~|~ s < N(\bi) \}$.

\begin{theorem} \label{degofgenerators}
Let $\Ix$ be the defining ideal of a set of $s$ points in $\pnk$
in generic position, and set 
$\mathcal{T} = \D \cup \left(\bigcup_{i=1}^{k} \D + e_i\right)$.
If $f$ is a generator of $\Ix$, then $\deg f \in \T$.  
In particular,
$
\Ix = \left\langle\bigoplus_{\bi \in \T} (\Ix)_{\bi} 
 \right\rangle.
$
\end{theorem}

\begin{proof}
The proof consists of two parts.

\noindent
{\it Step 1.}  We use Theorem \ref{generatorsprop} to show that
$\Ix = \left\langle\bigoplus_{\bi \in \T'} (\Ix)_{\bi} 
 \right\rangle$ where
\[\T' =  \bigcup_{\{l_1,\ldots,l_t\} \in \mathcal{P}(
[k])} \D + (e_{l_1} + \cdots + e_{l_t})\]
and $\mathcal{P}([k])$ denotes the power set of $[k]$.
It is enough to show that $\E \subseteq \T'$.  To do this, we need
to first show that $\D \subseteq \E$.  If $\bi \in \D$,
then $H_{\X}(\bi - 2e_l) < H_{\X}(\bi - e_l) \leq H_{\X}(\bi)= s$
for all $l\in [k]$.  So $\bi \in \D \subseteq \B$, but
$\bi \not\in \A$, hence $\bi \in \E$.

Suppose $\bj \in \E$.  Then there exists $\bi \in \D \subseteq \E$
such that $\bj \geq \bi$.  We can thus write
$\bj = (i_1+m_1,\ldots,i_k+m_k)$ where $\bi = (i_1,\ldots,i_k)$.  
If $m_1 = \cdots = m_k = 0$, then $\bj = \bi$ and hence
$\bj \in \E$.  So, suppose $m_l \geq 1$ for some $l \in [k]$.
If $m_l \geq 2$, then 
$H_{\X}(\bj) = H_{\X}(\bj - e_l) = H_{\X}(\bj-2e_l) = s $
since $\bj - 2e_l \geq \bi$.  But then $\bj \in \A$,
so $\bj \not\in \E$.  Hence, for each $l\in [k]$,
$m_l = 0$ or $1$.
So, if $m_{l_1} = \cdots = m_{l_t} = 1$,
and $0$ otherwise,  then $\bj = \bi + (e_{l_1}+\cdots + e_{l_t})
\in \T'$.

\noindent
{\it Step 2.} 
If $F$ is a generator of $\Ix$ with $\deg F = \bj$,
then the previous  step implies there exists $\bi \in \D$
such that $\bj = \bi + e_{l_1} + \cdots + e_{l_t}$ for
some subset $\{l_1,\ldots,l_t\} \subseteq [k]$.  We
wish to show that $t = 0$ or $1$, i.e., $\deg F = \bi$,
or $\deg F = \bi + e_l$ from some $l \in [k]$.

Let $\bi \in \D$ and let $\{l_1,\ldots,l_t\}$ be any subset of 
$[k]$ with $t \geq 2$.  Set $\bj = \bi + e_{l_1} + \cdots +
e_{l_t}$.  If we can show that
$(\Ix)_{\bj} = R_{e_{l_1}} (\Ix)_{\bj -e_{l_1}}
+ R_{e_{l_2}} (\Ix)_{\bj - e_{l_2}}$
then we shall be finished because this implies that
$(\Ix)_{\bj}$ contains no new generators.

By Remark \ref{nzd}, $\overline{x}_{l_2,0}$ is a 
non-zero divisor.  Set
$S = {\bf k}[x_{1,0},\ldots,\widehat{x}_{l_2,0},
\ldots,x_{k,n_k}] \cong 
R/(x_{l_2,0})$.  For each $\bt \in \N^k$ we have the
short exact sequence of vector spaces:
\[0 \longrightarrow (\Ix)_{\bt -e_{l_2}}
\stackrel{\times x_{l_2,0}}{\longrightarrow}
(\Ix)_{\bt} \longrightarrow 
\left({\Ix}/{x_{l_2,0}\Ix} \right)_{\bt} \cong
\left({(\Ix,x_{l_2,0})}/{(x_{l_2,0})}\right)_{\bt} 
\longrightarrow 0.\]
This gives $\dim_{\bf k} (\Ix)_{\bt} =
\dim_{\bf k} (\Ix)_{\bt -e_{l_2}} + \dim_{\bf k}
\left({(\Ix,x_{l_2,0})}/{(x_{l_2,0})}\right)_{\bt}$.

Since $\X$ is in generic position, 
$H_{\X}(\bj) = H_{\X}(\bj -e_{l_2}) = 
H_{\X}(\bj - e_{l_1}) = H_{\X}(\bj-e_{l_1}-e_{l_2}) = s$.
Thus, we can use the  short exact sequence
\[0 \longrightarrow R/\Ix(-e_{l_2})
\stackrel{\times \overline{x}_{l_2,0}}{\longrightarrow}
R/\Ix \longrightarrow R/(\Ix,x_{l_2,0}) \cong \frac{R/(x_{l_2,0})}
{(\Ix,x_{l_2,0})/(x_{l_2,0})} \longrightarrow 0\]
to show that $\left((\Ix,x_{l_2,0})/(x_{l_2,0})\right)_{\bj - e_{l_1}}
\cong S_{\bj - e_{l_1}}$ and 
$\left((\Ix,x_{l_2,0})/(x_{l_2,0})\right)_{\bj} \cong
S_{\bj}$.

Set $W = R_{e_{l_1}} (\Ix)_{\bj -e_{l_1}}
+ R_{e_{l_2}} (\Ix)_{\bj - e_{l_2}}$.  So $W \subseteq (\Ix)_{\bj}$.
Because
$x_{l_2,0}$ is a non-zero divisor
\[\dim_{\bf k} W = \dim_{\bf k} (\Ix)_{\bj -e_{l_2}} + \dim_{\bf k} W'
= \dim_{\bf k} (\Ix)_{\bj} - \dim_{\bf k} S_{\bj} + \dim_{\bf k} W'\]
where $W' =  \{ f''~|~ f = f'x_{l_2,0} + f'', ~f \in W \}$.
The vector space $W'$ can be viewed as a subspace of $S_{\bj}$.
It now suffices to show that 
$S_{\bj} \subseteq W'$
because then $\dim_{\bf k} W = \dim_{\bf k}(\Ix)_{\bj}$,
and thus, $W = (\Ix)_{\bj}$.  So, let $f \in S_{\bj}$.  Then
$f = f_0x_{l_1,0} + \cdots + f_{n_{l_1}}x_{l_1,n_{l_1}}$
with $f_i \in S_{\bj -e_{l_1}}$.  Since $S_{\bj -e_{l_1}}
\cong
\left((\Ix,x_{l_2,0})/(x_{l_2,0})\right)_{\bj - e_{l_1}}$, there
exists (with a slight abuse of notation)
$F_i \in (\Ix)_{\bj- e_{l_1}}$ such that $F_i
= g_ix_{l_2,0} + f_i$ .  But then 
$F_0x_{l_1,0} + \cdots + F_{n_{l_1}}x_{l_1,n_{l_1}} \in W$,
and hence, $f \in W'$.
\end{proof}

\begin{remark}
If $k = 1$ and $\X$ is a set of $s$ points in generic position,
then  we obtain the
well known result that 
$\Ix = \langle I_d\oplus I_{d+1} \rangle$ where
$d = \min\{i ~|~ \binom{n+ i}{i} > s\}.$  
\cite[Lemma 4.2]{GuMaRa2} is a proof of this theorem
in the special case
that $\X$ is a set of points
in generic position in $\popo$.  
If $\X$ is a set of points in generic position in $\pnk$,
and if we set $d_i := \min\left\{d ~\left|~ \binom{n_i + d}{d} \geq |\X|
\right\}\right.$ and
$D:=\max\{d_1+1,\ldots,d_k+1\}$, then the above
result implies that $\Ix$, considered as an $\N^1$-graded ideal
of $R$, is generated by forms of degree $\leq D$.  This is
extended in \cite{VT3} to show that $\operatorname{reg}(\Ix) =D$,
where $\operatorname{reg}(\Ix)$ is the Castelnuovo-Mumford
regularity of $\Ix$.
\end{remark}

\begin{corollary} \label{nomingen2}
Let $\X$ be a set of $s$ points in
generic position in $\pnk$ with 
Hilbert function $H_{\X}$.
If there exists $l,m \in [k]$ (not
necessarily distinct) and $\bj \in \N^k$ such
that $H_{\X}(\bj) = H_{\X}(\bj - e_l) = H_{\X}(\bj-e_l-e_m) = s$,
then $(\Ix)_{\bj}$ contains no generators of $\Ix$.
\end{corollary}

\begin{proof}
If $l = m$, then this is simply Theorem \ref{nomingen1}.
If $l \neq m$, then $H_{\X}(\bj - e_m) = s$ because
$\bj - e_l -e_m \leq \bj - e_m$ and $\X$ is in generic position.
Arguing as in Step 2 of Theorem \ref{degofgenerators}
we have $(\Ix)_{\bj} = R_{e_l}(\Ix)_{\bj-e_l} + 
R_{e_m} (\Ix)_{\bj -e_m}.$
\end{proof}

%%%%%%%%%%%%%%%%%%%%%%%%%%%%%%%%%%%%%%%%%%%%%%%%%%%%%%%%%%%%%%%%%%%%%%%

\section{On the value of $\nu(\Ix)$ for points in generic position}

In this section we study $\nu(\Ix)$, the minimal number of 
generators of $\Ix$, when $\X$ is a set
of points in generic position.  Unless specified otherwise, the
set of points under consideration will be non-degenerate, that is,
$|\X| > \max\{n_1,\ldots,n_k\}$.
We give an upper bound
on $\nu(\Ix)$ that can be calculated from $n_1,\ldots,n_k$ and $|\X| =s$.
We also show that calculating $\nu(\Ix)$ is equivalent
to calculating the dimensions of specific vector
spaces.  In some special cases, we  are able to compute
these dimensions.

So, suppose $\X$ is a non-degenerate set of $s$ points in generic position 
in $\pnk$.  By Theorem \ref{degofgenerators}, 
if   $\bi \in \D = \min\{ \bi \in \N^k ~|~ s < N(\bi)\}
= \min\{\bi ~|~ (\Ix)_{\bi} \neq 0\}$, by degree considerations
$(\Ix)_{\bi}$ cannot be generated by elements of smaller degree.
So the $\dim_{\bf k} (\Ix)_{\bi} = N(\bi) - H_{\X}(\bi)$ 
linearly independent elements of $(\Ix)_{\bi}$
must be generators of $\Ix$.  This gives a crude bound on $\nu(\Ix)$:
\[
\nu(\Ix) \geq \sum_{\bi \in \D} \dim_{\bf k} (\Ix)_{\bi} =
\sum_{\bi \in \D} \left(N(\bi) - s\right).
\]
By Theorem \ref{degofgenerators}, 
to compute $\nu(\Ix)$ it suffices to calculate
the number of generators of $\Ix$ in $(\Ix)_{\bj}$
for each $\bj \in \bigcup_{l=1}^k (\D + e_{l})$.

We wish to describe a subset of $\bigcup_{l=1}^k (\D + e_l)$
such that for each $\bj$ in this subset, $(\Ix)_{\bj}$ contains
no new generators of $\Ix$.  We introduce some 
suitable notation.  
For each $\bi \in \D$ set
\[\DD_{\bi}:= \{\bj \in \N^k ~|~ \bj \geq \bi\}\backslash
\{\bi,\bi+e_1,\bi+e_2,\ldots,\bi+e_k\}.\]
If follows that $\bj \in \DD_{\bi}$ if and only if 
$\bj - e_{l_1} - e_{l_2} \geq \bi$ for some not
necessarily distinct $l_1,l_2 \in [k]$.  

\begin{lemma} \label{subsetofj}
Let $\X$ be a non-degenerate set of $s$ points in $\pnk$.
With the notation as above, suppose
$\bj \in \left[\bigcup_{l=1}^k (\D + e_l)  \right] \cap 
\left[ \bigcup_{\bi \in \D} \DD_{\bi}\right]. $
Then $\Ix$ has no generator of degree $\bj$.
\end{lemma}

\begin{proof}
Since $\bj \in \DD_{\bi}$ for some $\bi \in \D$, $\bj -e_{l_1} - e_{l_2}
\geq \bi$ for some not necessarily distinct $l_1,l_2 \in [k]$.
Because $\X$ is in generic position, we have 
$H_{\X}(\bj - e_{l_1} - e_{l_2}) = H_{\X}(\bj - e_{l_1}) =
H_{\X}(\bj) = s$, and so the
conclusion follows from Corollary \ref{nomingen2}.
\end{proof}

Set 
\[\DD := 
\left[\bigcup_{l=1}^k (\D + e_l)  \right] \backslash 
\left[ \bigcup_{\bi \in \D} \DD_{\bi}\right].\]
Because of Lemma \ref{subsetofj}, to determine
$\nu(\Ix)$ it is enough to count the number of 
generators of $\Ix$ with degree $\bj \in \DD$. 

So, let $\bj \in \DD$.  
Since $\bj \in \bigcup_{l=1}^k (\D + e_{l})$, we can associate 
to $\bj$ a unique subset
$L_{\bj} := \{l_{1},\ldots,l_t\} \subseteq [k]$ such that 
$\bj \in \D + e_{l_m}$ for each $l_m \in L_{\bj}$ but
$\bj \not\in \D+e_l$ if $l \in [k] \backslash L_{\bj}$.  For
each $l_m \in L_{\bj}$ there then exists a unique $\bi_{l_m} \in \D$
such that $\bj = \bi_{l_m} + e_{l_m}$.  
So, for each $l_m \in L_{\bj}$ we can define 
$W_{l_m,\bi_{l_m}} := R_{e_{l_m}}(\Ix)_{\bi_{l_m}}
\subseteq (\Ix)_{\bj}.$
For each $\bj \in \DD$ we set
\[
W_{\bj} := W_{l_1,\bi_{l_1}} + \cdots + W_{l_t,\bi_{l_t}} 
= \sum_{l_m \in L_{\bj}} W_{l_m,\bi_{l_m}}
\subseteq
(\Ix)_{\bj}.
\]
Thus $W_{\bj}$ is the
subvector space of $(\Ix)_{\bj}$ that consists of all the 
forms in $\Ix$ of degree $\bj$ that come
from forms of lower degree in $\Ix$.
The number
of new generators of $\Ix$ of degree $\bj$  with $\bj \in \DD$ is 
then
\[
\dim_{\bf k}(\Ix)_{\bj} - \dim_{\bf k} 
(W_{l_1,\bi_{l_1}} + \cdots + W_{l_t,\bi_{l_t}}) = N(\bj) -s - 
\dim_{\bf k} W_{\bj}.
\]
We summarize this discussion with the following theorem.

\begin{theorem} \label{numberofgenerators}
Let $\X$ be a non-degenerate set of $s$ points in generic position 
in $\pnk$.  With 
the notation as above
\begin{eqnarray*}
\nu(\Ix) & = & \sum_{\bi \in \D} \dim_{\bf k} (\Ix)_{\bi}
+ \sum_{\bj \in \DD} \left(\dim_{\bf k} (\Ix)_{\bj}
- \dim_{\bf k} W_{\bj} \right) \\
& = &
\sum_{\bi \in \D} \left(N(\bi)-s \right) + 
\sum_{\bj \in \DD}
\left(N(\bj) -s - \dim_{\bf k} W_{\bj}\right).
\end{eqnarray*}
\end{theorem}

Computing $\nu(\Ix)$ is thus equivalent to computing
$\dim_{\bf k} W_{\bj}$ for each $\bj \in \DD$.  
Arguing as in \cite[Proposition 7]{GO2} one has the following
lower bounds:

\begin{lemma} \label{simplebounds}
Suppose $\bj \in \DD$ and $L_{\bj} =
\{l_1,\ldots,l_t\}$.  For every $l_m \in L_{\bj}$
\[
\dim_{\bf k} W_{\bj} \geq \dim_{\bf k} W_{l_m},\bi_{l_m} \geq
2 \dim_{\bf k} (\Ix)_{\bi_{l_m}}.
\]
\end{lemma}

Combining Lemma \ref{simplebounds} 
with Theorem \ref{numberofgenerators} gives us
an upper bound on $\nu(\Ix)$.

\begin{corollary}
Let $\X$ be a non-degenerate set of $s$ points in
generic position in $\pnk$.  With
the notation as above
\begin{eqnarray*}
\nu(\Ix) & \leq &\sum_{\bi \in \D} \left(N(\bi) - s\right) +
\sum_{\bj \in \DD}
\left(N(\bj) - 2N(\bj -e_{l_1}) + s\right)
\end{eqnarray*}
where $l_1 \in L_{\bj} = \{l_1,\ldots,l_t\}$ for $\bj \in \DD$.
\end{corollary}

\begin{proof}
For each $\bj \in \DD$, let $l_1 \in L_{\bj}$.
Then $\dim_{\bf k} W_{\bj} \geq 2\dim_{\bf k} (\Ix)_{\bj-e_{l_1}}
= 2N(\bj-e_{l_1}) - 2s.$  Now use Theorem \ref{numberofgenerators}.
\end{proof}

In general, calculating $\dim_{\bf k} W_{\bj}$ for each $\bj \in \DD$
is a very difficult problem.  Indeed, if $L_{\bj} = \{l_1,\ldots,l_t\}$,
then there is no {\it a priori} formula for calculating 
$\dim_{\bf k} W_{l_m,\bi_{l_m}} = 
\dim_{\bf k} (R_{e_{l_m}}(\Ix)_{\bi_{l_m}})$ for each $l_m \in L_{\bj}$.  
The problem is further complicated when 
 $|L_{\bj}| = t \geq 2$ because then 
we need to know how $W_{l_m,{\bi}_{l_m}}$ and $W_{l_n,{\bi}_{l_n}}$ 
intersect in $(\Ix)_{\bj}$
for each $l_n, l_m \in L_{\bj}$.

However, under some extra hypotheses on either $s = |\X|$ or
$n_1,\ldots,n_k$ we can be quite explicit about $\dim_{\bf k} W_{\bj}$
for {\it some} $\bj \in \DD$.  The remaining results of this
section are of this vein.  

\begin{lemma}	\label{growthinp1}
Let $\X \subseteq \pr^1 \times \cdots \times \pr^1$ be any
finite set of points.  Then, for any $\bi \in \N^k$ and $l \in [k]$,
\[
\dim_{\bf k} \left( R_{e_l}(\Ix)_{\bi} \right)
= 2 \dim_{\bf k} (\Ix)_{\bi} - \dim_{\bf k} (\Ix)_{\bi -e_l}
\]
where $\dim_{\bf k} (\Ix)_{\bi -e_l} = 0$ if $\bi -e_l \not\geq 0$.
\end{lemma}

\begin{proof}
The proof for the case $\X \subseteq \popo$ given in
\cite[Lemma 2.3]{GuMaRa1} 
can be extended to this case.
\end{proof}

\begin{theorem} \label{growth1}
Let $\X \subseteq \pr^1 \times \cdots \times \pr^1$ be a set
of $s>1$ points in generic position.  With the notation as above,
suppose $\bj \in \DD$ with $L_{\bj} = \{l\}$.  Then
\[
\dim_{\bf k} W_{\bj} =
\left\{
\begin{array}{ll}
\dim_{\bf k} (\Ix)_{\bj} = N(\bj) - s & \mbox{if $N(\bj -2e_l) = s$.} \\
2\dim_{\bf k} (\Ix)_{\bj-e_l} = 2N(\bj-e_l) - 2s & 
 \mbox{if $N(\bj -2e_l) < s$.}
\end{array}
\right.
\]
\end{theorem}

\begin{proof}
The hypothesis $L_{\bj} = \{l\}$ implies $\bj - e_{l} = \bi \in \D$
but $\bj - e_m \not\in \D$ for any $m \in [k]\backslash\{l\}$.
Hence $\dim_{\bf k} W_{\bj} = \dim_{\bf k} (R_{e_l}
(\Ix)_{\bi})$.  Since $\bi \in \D$, there does not exists an $\bi' \in \D$
such that $\bi -e_l \geq \bi'$, and hence
$\dim_{\bf k} (\Ix)_{\bi -e_l} = 0$.  By Lemma \ref{growthinp1}
we thus have
$\dim_{\bf k} W_{\bj} = 2\dim_{\bf k} (\Ix)_{\bi} = 2N(\bj-e_l)-2s.
$
The reader can now verify that $2N(\bj-e_l) -2s = N(\bj) -s$ if
$s = N(\bj -2e_l)$.
\end{proof}

\begin{theorem} 	\label{growth2}
Let $\X$ be a non-degenerate set of $s$ points in generic position
in $\pnk$
and $\bj \in \DD$.
\begin{enumerate}
\item[$(i)$]  If there exists $l \in [k]$ such
that $N(\bj - 2e_l) = s$, then $\dim_{\bf k} W_{\bj} =
\dim_{\bf k} (\Ix)_{\bj}$.
\item[$(ii)$]  If $L_{\bj} = \{l\}$ and $s = N(\bj -e_l) -1$,
then $\dim_{\bf k} W_{\bj} = n_l + 1$.
\end{enumerate}
\end{theorem}

\begin{proof} $(i)$  Since $N(\bj - 2e_l) = s$, $H_{\X}(\bj -2e_l)
= H_{\X}(\bj -e_l) = s$.  Now apply Theorem \ref{generators}.

$(ii)$  We are given that $\bj - e_l \in \D$ and 
$\dim_{\bf k} (\Ix)_{\bj -e_l} = N(\bj-e_l) -s =1$.  So let $F$
be a basis for $(\Ix)_{\bj -e_l}$.  Then
$x_{l,0}F,\ldots,x_{l,n_l}F$ form a linearly independent
basis of $W_{\bj} = R_{e_l} (\Ix)_{\bj-e_l}$.
\end{proof}
%%%%%%%%%%%%%%%%%%%%%%%%%%%%%%%%%%%%%%%%%%%%%%%%%%%%%%%%%%%%%%%%%%%%%%%

\section{On the expected value of $\nu(\Ix)$}

Let $\X$ be a non-degenerate set of points in generic position
in $\pnk$.  In this section we are interested in determining
if there is an expected value for $\nu(\Ix)$.  After showing
that $\nu(\Ix)$ is constant on some open subset of $(\pnk)^s$,
we give a lower bound for this value.  When $k = 1$
the resulting lower bound is conjectured to equal $\nu(\Ix)$ on
some non-empty open subset of $(\pr^n)^s$ by the
Ideal Generation Conjecture.  Therefore,
it seems natural to expect that our generalized
lower bound equals $\nu(\Ix)$ on some non-empty open
subset of $(\pnk)^s$, thus generalizing the IGC to points in $\pnk$.  However,
although we have found computationally that in many cases $\nu(\Ix)$
equals the lower bound, 
we show that there exist $s$ and $n_1,\ldots,n_k$ for which
$\nu(\Ix)$ is always larger than this bound.
We continue to use the notation of the previous sections.

If $(P_1,\ldots,P_s) \in (\pnk)^s$ is such that $P_1,\ldots,P_s$
are distinct points, then we shall write $I(P_1,\ldots,P_s)$
to denote the defining ideal of $\{P_1,\ldots,P_s\} \subseteq
\pnk$.  Furthermore, if $\bj \in \D$ with $L_{\bj} = 
\{l_1,\ldots,l_t\}$, then we write $W(P_1,\ldots,P_s)_{\bj}$
for the vector space
$W(P_1,\ldots,P_s)_{\bj} := R_{e_{l_1}}
I_{\bj-e_{l_1}} + \cdots + R_{e_{l_t}} I_{\bj -e_{l_t}}
\subseteq I_{\bj}$
where $I = I(P_1,\ldots,P_s)$.  

\begin{theorem}
 Let $s > \max\{n_1,\ldots,n_k\}$.
Then there exists an open set $U \subseteq (\pnk)^s$
such that if $(P_1,\ldots,P_s) \in U$, then
$\dim_{\bf k} W(P_1,\ldots,P_s)_{\bj}$ is the maximum
possible for all $\bj \in \DD$.  In particular, 
$\nu(I(P_1,\ldots,P_s))$ is constant for
all $(P_1,\ldots,P_s) \in U$.
\end{theorem}

\begin{proof}
It is enough to show that for each $\bj \in \DD$, there
exists an open subset $U_{\bj} \subseteq (\pnk)^s$
with the property that if $(P_1,\ldots,P_s) \in U_{\bj}$,
then $\dim_{\bf k} W(P_1,\ldots,P_s)_{\bj}$ is maximal.  Then,
since $|\DD| < \infty$, the desired open set is
$U = \bigcap_{\bj \in \DD} U_{\bj}$.

So, let $\bj \in \DD$ and suppose $L_{\bj} = \{l_1,\ldots,l_t\}$.
For each $l_m \in L_{\bj}$ set $\bi_{l_m} := \bj - e_{l_m}$.  Let
$W \subseteq (\pnk)^s$ denote the open set of Theorem \ref{genericpoints}
consisting of the $s$ distinct points in generic position.  Then,
using the proof of the claim found after Remark 2.8 in \cite{GM},
we can show that there exists an open set $U_{l_m} \subseteq W$
such that for all $(P_1,\ldots,P_s) \in U_{l_m}$,
$\dim_{\bf k} R_{e_{l_m}} 
I(P_1,\ldots,P_s)_{\bi_{l_m}}$
is the maximum possible.

If we let $G_{l_m,1},\ldots,G_{l_m,N(\bi_{l_m})-s}$ denote the $N(\bi_{l_m})
-s = \dim_{\bf k} I(P_1,\ldots,P_s)_{\bi_{l_m}}$ distinct basis elements
of $I(P_1,\ldots,P_s)_{\bi_{l_m}}$, then the elements
\[
\left\{x_{l_m,i}G_{l_m,j} ~\left|~ 0 \leq i \leq n_{l_m}, 1 \leq j \leq 
N(\bi_{l_m}) -s \right\}\right.\]
 generate 
$R_{e_{l_m}} I(P_1,\ldots,P_s)_{\bi_{l_m}}$.  
Set $M_{l_m} = (n_{l_m}+1)(N(\bi_{l_m})-s)$
and form the $M_{l_m} \times N(\bj)$ matrix $\mathcal{M}_{l_m}$
which expresses how the $x_{l_m,i}G_{l_m,j}$'s are linear combinations
of the $N(\bj)$ monomials of degree $\bj$.  Since rank $\mathcal{M}_{l_m}
= \dim_{\bf k} 
R_{e_{l_m}}I(P_1,\ldots,P_s)_{\bi_{l_m}}$,
this rank is maximal for all $(P_1,\ldots,P_s) \in U_{l_m}$.

Let $\mathcal{M}$ be the $\left(\sum_{l_m \in L_{\bj}} M_{l_m}\right)
\times N(\bj)$ matrix
$\mathcal{M} := \bmatrix
\mathcal{M}_{l_1}\\
\vdots \\
\mathcal{M}_{l_t}
\endbmatrix
.$
Then the rank of $\mathcal{M}$ is equal to $\dim_{\bf k} 
W(P_1,\ldots,P_s)_{\bj}$.  The rank of $\mathcal{M}$ will therefore
assume its maximal value on some open subset $U_{\bj} \subseteq
\bigcap_{l_m \in L_{\bj}} U_{l_m}$.  This is the desired set $U_{\bj}$.
\end{proof}

We can give a lower bound on $\nu(\Ix)$ by using Theorem 
\ref{numberofgenerators} and bounds on  $\dim_{\bf k} W_{\bj}$
for each $\bj \in \DD$.  For
each $l_m \in L_{\bj}$ the dimension of $W_{l_m,\bi_{l_m}}$
is bounded by
\begin{equation}\label{bound1}
\dim_{\bf k}  W_{l_m,{\bi}_{l_m}}  \leq
\min \left\{\dim_{\bf k} (\Ix)_{\bj}, ~
(n_{l_m}+1)\dim_{\bf k}(\Ix)_{\bi_{l_m}}\right\}.
\end{equation}
Furthermore, if $l_m,l_n \in L_{\bj}$ and $l_m \neq l_n$,
then 
\[\dim_{\bf k} (W_{l_m,\bi_{l_m}} + W_{l_n,\bi_{l_n}}) \leq
\dim_{\bf k} W_{l_m,\bi_{l_m}} + \dim_{\bf k} W_{l_n,\bi_{l_n}}.\]
We thus arrive at the following upper bound for $W_{\bj}$:
\begin{equation}\label{bound2}
\dim_{\bf k} W_{\bj} \leq \min\left\{\dim_{\bf k} 
(\Ix)_{\bj}, ~\sum_{l_m \in L_{\bj}}
(n_{l_m}+1) \dim_{\bf k} (\Ix)_{\bi_{l_m}}\right\}.
\end{equation}
Since the values of $\dim_{\bf k} (\Ix)_{\bj}$
and $\dim_{\bf k} (\Ix)_{\bi_{l_m}}$ are known because $\X$
is in generic position, combing the 
above upper bound with Proposition \ref{numberofgenerators}
results in the following lower bound: 
\begin{theorem}
Let $\X$ be a non-degenerate set of $s$ points in generic
position in $\pnk$, and set
{\small
\[v(s;n_1,\ldots,n_k):=\sum_{\bi \in \D} (N(\bi) -s)
+ \sum_{\bj  \in \DD}
\max\left\{0, ~N(\bj)-s - \sum_{l_m \in L_{\bj}} 
(n_{l_m}+1)(N(\bi_{l_m}) - s)\right\}. 
\]}Then $\nu(\Ix) \geq v = v(s;n_1,\ldots,n_k)$.
\end{theorem}

\begin{remark}
If $k=1$, then 
\[v = v(s;n) = \binom{d+n}{n}-s+\max\left\{0,\binom{d+1+n}{n}-s-(n+1)
\left(\binom{d+n}{n}-s\right)\right\}\] 
where $d = \min\left\{i ~\left| \binom{n+i}{i} > s\right\}\right.$.
The Ideal Generation Conjecture conjectures that $\nu(\Ix) = v$
on some non-empty open subset of $(\pr^n)^s$.  Although
known to be true in some cases (for $n=2$ see \cite{GM},
for $n=3$ see \cite{B}, and for $s \gg n$ see
\cite{HS}) the conjecture remains open in general.  The conjecture
was formulated using the heuristic argument that ``generically''
$\dim_{\bf k} W_{\bj}$ should be as large as possible, thus implying
equality in the bounds $(\ref{bound1})$ and $(\ref{bound2})$.
It seems natural to extend this heuristic argument to points
in $\pnk$ to generalize the Ideal Generation Conjecture by expecting
that $\nu(\Ix) =v(s;n_1,\ldots,n_k)$ 
for some non-empty open set of $(\pnk)^s$.  But
as we show at the end of this section, sometimes $\nu(\Ix) > v$ 
if $k \geq 2$.
\end{remark}

In \cite{GuMaRa2,GuMaRa3} Giuffrida, et al. computed the minimal
free resolution of points in generic position in  $\popo$, and 
in particular, showed that $\nu(\Ix) = v$ on some non-empty open
set of $(\popo)^s$.  Since their notation
and approach is different than ours, for the convenience of the reader
we make this connection more transparent.

\begin{theorem}
There exists a non-empty open subset $U \subseteq (\popo)^s$ with $s \geq 2$
such that for all $(P_1,\ldots,P_s) \in U$, the points $\{P_1,\ldots,P_s\}$
are in generic position, and $\nu(I(P_1,\ldots,P_s))=v(s;1,1)$.
\end{theorem}

\begin{proof}
It suffices to show that the bound (\ref{bound2})
for $\dim_{\bf k} W_{\bj}$ is in fact an equality for each $\bj \in \DD$.
If $\bj = (i,j) \in \DD$,  there are two possibilities: $|L_{\bj}| = 1$
or $2$.  
In the former, by Theorem \ref{growth1}
we have equality in $(\ref{bound2})$.

For the second case, by \cite[Theorem 4.3]{GuMaRa3} there
exists a non-empty subset $U \subseteq (\popo)^s$ such that for all
$\bj = (i,j) \in \DD$ with $|L_{\bj}| =2$
and for each $(P_1,\ldots,P_s) \in U$, we have
\[\dim_{\bf k} (I(P_1,\ldots,P_s))_{i,j} - 
\dim_{\bf k} W(P_1,\ldots,P_s)_{i,j} = \max\{0,-d_{i,j}\}.\]
Here, $d_{i,j}$ is the $(i,j)${th} entry of what \cite{GuMaRa3}
call the second difference Hilbert matrix of $\X = \{P_1,\ldots,P_s\}$ 
which is computed from the Hilbert function on 
$\X$.  Since $|L_{\bj}| = 2$ and because
$\X$ is in generic position, $H_{\X}$,
written as a matrix, has the form
\[
\begin{tabular}{r|cccc}
 & & $_{(j-2)}$ & $_{(j-1)}$ & $_{j}$ \\
\hline
 & & $\vdots$ & $\vdots$ & $\vdots$ \\
$_{(i-2)}$ & $\cdots$ & $(i-1)(j-1)$ & $(i-1)j$ & $(i-1)(j+1)$\\
$_{(i-1)}$ & $\cdots$ & $i(j-1)$ & $ij$ & $s$ \\
$_{i}$ & $\cdots$ & $(i+1)(j-1)$ & $s$ & $s$
\end{tabular}
\]
This local description of the Hilbert function, and the 
definition of $d_{i,j}$ on page 422 of \cite{GuMaRa3} gives
\begin{eqnarray*}
-d_{i,j} & =& \left[(i+1)(j+1)-s\right] -
\left[2((i+1)j -s) + 2((j+1)i-s)\right] \\ 
& = & \dim_{\bf k} (\Ix)_{i,j} - 2\dim_{\bf k} (\Ix)_{i-1,j} - 2\dim_{\bf k}
(\Ix)_{i,j-1}.
\end{eqnarray*}
We thus have equality in $(\ref{bound2})$.
\end{proof}

We now show that $\nu(\Ix)$ may not equal $v =
v(s;n_1,\ldots,n_k)$ in general.
We begin by showing that any example
of points $\X$ with $\nu(\Ix) >v$
can be extended to an infinite family of examples.

\begin{lemma} \label{makecounterexamples}
Suppose  that for every non-degenerate set 
$\X$ of $s$ points  in generic position
in $\pnk$ we have $\nu(I_{\X}) > v(s;n_1,\ldots,n_k)$.
If $\X'$ is any non-degenerate set of s points in generic
position in $\pnk \times \pr^{m_1} \times \cdots \times
\pr^{m_l}$, then
$\nu(I_{\X'}) > v(s;n_1,\ldots,n_k,m_1,\ldots,m_l)$.
\end{lemma}

\begin{proof}
Let $\X'$ be a set of $s$ points in generic position
in $\pnk \times \pr^{m_1} \times \cdots \times \pr^{m_l}$.
Let $I_{\X'}$ be the associated ideal and set
$
I := \bigoplus_{(i_1,\ldots,i_k) \in \N^k} 
(I_{\X'})_{(i_1,\ldots,i_k,0,\ldots,0)}.
$
Then $I$ is isomorphic to an ideal $I_{\X} \subseteq
{\bf k}[x_{1,0},\ldots,x_{1,n_1},\ldots,x_{k,0},\ldots,x_{k,n_k}]$
where $I_{\X}$ is the defining ideal of a set $\X$ of $s$
points in generic position in $\pnk$.

By hypothesis,
there exists $\bj = (j_1,\ldots,j_k) \in \N^k$ such
that $(I_{\X})_{\bj}$ contains a generator that has not been accounted
for by $v(s;n_1,\ldots,n_k)$.  Hence 
$(I_{\X'})_{(j_1,\ldots,j_k,0,\ldots,0)} \cong (\Ix)_{\bj}$
contains a generator of $I_{\X'}$ that is not expected,
and thus $\nu(I_{\X'})$ will be strictly 
larger than $v(s;n_1,\ldots,n_k,m_1,\ldots,m_l)$.
\end{proof}

We now give a case where $\nu(\Ix)$ fails to agree with the lower bound.

\begin{theorem} \label{badexamples}
Let $\X$ be three points in generic position in $\pr^1 \times \pr^1
\times \pr^1$.  Then  $\nu(\Ix) > v(3;1,1,1)$.
\end{theorem}

\begin{proof}
It is enough to show the existence of some $\bj \in \DD$ 
for which we have a strict inequality in $(\ref{bound2})$
for $\dim_{\bf k} W_{\bj}$.
Now $\{ (1,1,0),(1,0,1),(0,1,1)\} \subseteq \D = 
\min\{\bi\in\N^3 ~|~N(\bi) > 3\}$, 
and so $\bj := (1,1,1) \in \DD$
with $L_{\bj} = \{1,2,3\}$.  The expected dimension of
$W_{\bj}$ is 
\[\min\left\{\dim_{\bf k} (\Ix)_{\bj},
\sum_{i\in\{1,2,3\}} 2\dim_{\bf k} (\Ix)_{\bj -e_i}\right\} = 
\min\left \{ 8 -3, \sum_{i\in\{1,2,3\}} 2(4-3)\right\} = 5,\]
or equivalently, we expect $(\Ix)_{\bj}$ to contain no generators.

However, we claim $\dim_{\bf k} W_{\bj} \leq 4$, and hence,
$(\Ix)_{\bj}$ contains a new generator.  Let $P_1,P_2,P_3$
be the distinct points of $\X$, and after a linear change of variables
in each set of coordinates, we can assume $
P_1  =  [1:0] \times [1:0] \times [1:0]$, 
$P_2   =  [1:a_1] \times [1:a_2] \times [1:a_3]$,
and $P_3   =  [1:b_1] \times [1:b_2] \times [1:b_3]$
with $a_i \neq b_i$ for $i =1,2,3$.  Because $\X$
is in generic position $\dim_{\bf k} (\Ix)_{1,1,0} 
=\dim_{\bf k} (\Ix)_{1,0,1} = \dim_{\bf k} (\Ix)_{0,1,1} = 1$.
To find a basis for each of these vector spaces, it suffices to find a form
of the proper degree in $\Ix$.  From our description of the points
we can find such forms:
\begin{eqnarray*}
F_1 & := & (a_2b_1 - a_1b_2)x_1y_1 + a_2b_2(a_1 - b_1)x_1y_0 + 
a_1b_1(b_2-a_2)x_0y_1 \in (\Ix)_{1,1,0} \\
F_2 & := & (a_3b_1 - a_1b_3)x_1z_1 + a_3b_3(a_1 - b_1)x_1z_0 + 
a_1b_1(b_3-a_3)x_0z_1 \in (\Ix)_{1,0,1}\\
F_3 & := & (a_2b_3 - a_3b_2)y_1z_1 + a_3b_3(b_2 - a_2)y_1z_0 + 
a_2b_2(a_3-b_3)y_0z_1 \in (\Ix)_{0,1,1}.
\end{eqnarray*}
It  follows that $z_0F_1,z_1F_1,y_0F_2,y_1F_2,x_0F_3,x_1F_3$
generate the vector space $W_{\bj}$.  A routine calculation
will now verify that 
\begin{eqnarray*}
 a_1b_1(a_1-b_1)x_0F_3 & = & 
(a_1-b_1)[a_3b_3z_0F_1 -a_2b_2y_0F_2]+ (a_3b_1-a_1b_3)z_1F_1\\
&&
- (a_2b_1 - b_2a_1)y_1F_2 \\
 (a_1-b_1)x_1F_3 & = & (b_2-a_2)y_1F_2 + (a_3-b_3)z_1F_1. 
\end{eqnarray*}
Thus, $x_0F_3,x_1F_3$ are in the vector space spanned by
$z_0F_1,z_1F_1,y_0F_2,y_1F_2$, whence $\dim_{\bf k} W_{\bj}
\leq 4 < 5 = $ the expected dimension.
\end{proof}

With this result we can construct examples with $\nu(\Ix)$ arbitrarily larger
than $v(s;n_1,\ldots,n_k)$.

\begin{corollary}
Let $\X$ be three points in generic position in $\pr^1 \times \cdots
\times \pr^1$ ($k \geq 3$ times).  Then 
$\nu(\Ix) \geq v(3;1,\ldots,1) + \binom{k}{3}$.
\end{corollary}

\begin{proof}
There are $\binom{k}{3}$ tuples $\bi \in \N^k$ which have exactly three $1$'s
and $k-3$ zeroes.  Let $\bi$ be such a tuple, and suppose that the three 
$1$'s are in $i_1$th, $i_2$th,
and $i_3$th position.  If $\pi_{i_1,i_2,i_3}:\pr^1
\times \cdots \times \pr^1 \rightarrow \popo
\times \pr^1$ is the projection map onto the  $i_1$th, $i_2$th,
and $i_3$th coordinates, then $\Y = \pi_{i_1,i_2,i_3}(\X) \subseteq
\popo \times \pr^1$  is
a set of three points in generic position. 
Hence, $(\Ix)_{\bi} \cong (I_{\Y})_{1,1,1}$.  But by Theorem
\ref{badexamples} $(I_{\Y})_{1,1,1}$ contains at least one
generator not accounted for by $v(3;1,1,1)$, and thus, 
$(\Ix)_{\bi}$ has a generator not counted by $v(3;1,\ldots,1)$.
\end{proof}

Using CoCoA \cite{C} we have computed $\nu(\Ix)$ in the 
following ranges:
\[
\begin{array}{rll}
k = 2 & 1 \leq n_1 \leq n_2 \leq 5 & n_2 < s \leq 20 \\
k = 3 & 1 \leq n_1 \leq n_2 \leq n_3 \leq  5 & n_3 < s \leq 10 \\
k = 4 & 1 \leq n_1 \leq n_2 \leq n_3 \leq n_4 \leq 5 & n_4
< s \leq 10.
\end{array}
\]
Besides the example of Theorem \ref{badexamples} 
(and those examples that are a consequence of Lemma \ref{makecounterexamples})
we found that 
\[\nu(\Ix) > v(1+n+n;1+n+n) ~~\mbox{for $1 \leq n \leq 7$}.\]
From this data it appears that $\nu(\Ix) > v(1+n+n;1,n,n)$
for all $n$.
Notice that the example of Theorem \ref{badexamples}
is also part of this family.  Using CoCoA we found 
that in each of these cases  $\dim_{\bf k} W_{1,1,1}$ is
smaller than the expected dimension.  

We point out, however, that in every other case  the computed
value of $\nu(\Ix)$ agrees  with $v(s;n_1,\ldots,n_k)$.
These computations leads us
to believe that $\nu(\Ix) = v(s;n_1,\ldots,n_k)$ 
for a large number  $s$ and $n_1,\ldots,n_k$.
Moreover, we know of no counterexamples when $k \leq 2$.
We conclude by giving some questions inspired by our computer examples.

\begin{question}
For $s = (1+n+n)$  points in generic position in 
$\pr^1 \times \pr^n \times \pr^n$ is $\nu(\Ix)$ always
the larger $v(s;1,n,n)$?
Is this family of examples the only 
family where the lower bound  fails to hold?  If not,
can we classify all $s$ and $n_1,\ldots,n_k$ for which $\nu(\Ix)
\neq v(s;n_1,\ldots,n_k)$?  Does the lower bound 
value always hold in the case $k \leq2$?  How should a generalized
Ideal Generalization Conjecture be formulated to account
for these examples?
\end{question}

%%%%%%%%%%%%%%%%%%%%%%%%%%%%%%%%%%%%%%%%%%%%%%%%%%%%%%%%%%%%%%%%%%%%%%%

\end{document}